\documentclass[12pt]{amsart}
\input{epsf}
\usepackage{amssymb,latexsym}
\advance\textheight by \topskip
\numberwithin{equation}{section}
\newtheorem{theorem}{Theorem}[section]

\newtheorem{example}[theorem]{Example}

\begin{document}

\def\sof{\hfill\rule{2mm}{2mm}}
\def\ls{\leqslant}
\def\gs{\geqslant}
\def\SS{\frak S}
\def\T{\mathcal{T}}
\def\qq{{\bold q}}
\def\txx{{\tfrac1{2\sqrt{x}}}}

\title[Restricted permutations and Chebyshev polynomials]
{{\bf Restricted permutations and Chebyshev polynomials}}

\author[T. Mansour and A. Vainshtein]
{T. Mansour$^\ast$ and A. Vainshtein$^{\dag}$}
\maketitle
\begin{center}
$^\ast$ LABRI, Universit\'e Bordeaux\\
       351 cours de la Lib\'eration, 33405 Talence Cedex,
       France\\[4pt]

$^\dag$ Department of Mathematics and Department of Computer Science\\
University of Haifa, Haifa, Israel 31905\\[4pt]
{\tt toufik@labri.fr}, {\tt alek@mathcs.haifa.ac.il}
\end{center}

%===========================================================================

\begin{abstract}  We study generating functions for the number of permutations in $\SS_n$
subject to two restrictions. One of the restrictions belongs to $\SS_3$,
while the other belongs to $\SS_k$. It turns out that in a large
variety of cases
the answer can be expressed via Chebyshev polynomials of the second kind.

\end{abstract}
\ \\

{\Large\bf Contents:}\ \\

\contentsline {chapter}{\tocchapter {}{1}{{\bf
Introduction}\dotfill}}{2} \contentsline {section}{\tocsection
{}{1.1}{{\bf Pattern avoidance}\dotfill}}{2} \contentsline
{section}{\tocsection {}{1.2}{{\bf Chebyshev polynomials of the
second kind}\dotfill}}{2} \contentsline {section}{\tocsection
{}{1.3}{{\bf Preliminaries}\dotfill}}{3}

\contentsline {chapter}{\tocchapter {}{2}{{\bf Transfer
matrices}\dotfill}}{3} \contentsline {section}{\tocsection
{}{2.1}{{\bf Generating trees}\dotfill}}{4} \contentsline
{section}{\tocsection {}{2.2}{{\bf Dyck paths}\dotfill}}{5}

\contentsline {chapter}{\tocchapter {}{3}{{\bf Continued fractions
}\dotfill}}{6}

\contentsline {chapter}{\tocchapter {}{4} {\bf
                   Block decompositions\dotfill}}{7}
\contentsline {chapter}{\tocchapter {}{5} {\bf
                   Counting $132$-restricted permutations \dotfill}}{10}
\contentsline{section}{\tocsection{}{5.1}{{\bf Avoiding $132$ and
another pattern}\dotfill}}{10}
\contentsline{section}{\tocsection{}{5.2}{{\bf Avoiding $132$ and
containing another pattern exactly once}\dotfill}}{10}
\contentsline{section}{\tocsection{}{5.3}{{\bf Containing $132$
exactly once and avoiding another pattern}\dotfill}}{12}
\contentsline{section}{\tocsection{}{5.4}{{\bf Containing $132$
and another pattern exactly once}\dotfill}}{13}
\contentsline{section}{\tocsection{}{5.5}{{\bf
Generalizations}\dotfill}}{14}

\contentsline{chapter}{\tocchapter{}{6}{\bf Counting
$321$-restricted permutations\dotfill}}{16}

%==========================================================================
\setcounter{section}{1} \setcounter{subsection}{0}
\setcounter{theorem}{0}
\section*{\bf 1. Introduction}

\subsection*{1.1. Pattern avoidance}
Let $\pi=(\pi_1,\dots,\pi_n)\in\SS_n$ and $\tau\in\SS_k$ be two
permutations. An {\it
occurrence\/} of $\tau$ in $\pi$ is a subsequence
$\pi'=(\pi_{i_1}, \dots,\pi_{i_k})$ such that
$1\ls i_1<i_2<\dots<i_k\ls n$ and $\pi'$
is order-isomorphic to $\tau$; in such a context $\tau$ is usually
called a {\it pattern\/}. We say that $\pi$ {\it avoids\/} $\tau$,
or is $\tau$-{\it avoiding\/}, if there is no occurrence of $\tau$
in $\pi$. For example, let $\pi=83176254$, and
$\tau_1= 1234$, $\tau_2= 1243$, $\tau_3= 4213$.
Then it is easy to see that $\pi$ avoids $\tau_1$, contains exactly one
occurrence of $\tau_2$, namely, $1254$, and contains six occurrences
of $\tau_3$, which are $8317$,  $8316$, $8315$,
$8314$, $8325$, $8324$.

The set of all $\tau$-avoiding permutations in $\SS_n$
is denoted $\SS_n(\tau)$. For an arbitrary finite collection of
patterns $T$, we say that $\pi$ avoids $T$ if $\pi$ avoids any
$\tau\in T$; the corresponding subset of $\SS_n$ is denoted
$\SS_n(T)$. By $F_\tau$ and $F_T$ we denote the corresponding
(ordinary) generating functions
$F_\tau(x)=\sum_{n=0}^\infty|\SS_n(\tau)|x^n$ and
$F_T(x)=\sum_{n=0}^\infty|\SS_n(T)|x^n$.

The following simple symmetry arguments allow to decrease the
variety of different $T$'s to be considered. Define the {\it reversal\/}
and the {\it complementation\/} operations as follows:
$r(\pi_1,\pi_2,\dots,\pi_n)=(\pi_n,\pi_{n-1},\dots,\pi_1)$,
$c(\pi_1,\pi_2,\dots,\pi_n)=(n+1-\pi_1,n+1-\pi_2,\dots,n+1-\pi_n)$.
It is easy to see that the following four statements are equivalent:

(i) $\pi$ avoids $\tau$;

(ii) $r(\pi)$ avoids $r(\tau)$;

(iii) $c(\pi)$ avoids $c(\tau)$;

(iv) $\pi^{-1}$ avoids $\tau^{-1}$.

\noindent More generally, denote by $G$ the group of transformations generated
by $r$, $c$, and the usual group inverse operation (it is easy to see
that $G$ is isomorphic to the dihedral group $D_8$). The for any
$g\in G$ one has $F_T(x)=F_{g(T)}(x)$.

The first paper devoted
entirely to the study of permutations avoiding certain patterns
({\it restricted permutations\/}) appeared in 1985 (see
\cite{SS}). Currently there exist more than fifty papers on this
subject.

\subsection*{1.2. Chebyshev polynomials of the second kind}
{\it Chebyshev polynomials of the second kind\/} (in what follows
just Chebyshev polynomials) are defined by
$$
U_r(\cos\theta)=\frac{\sin(r+1)\theta}{\sin\theta}
$$
for $r\gs0$. Evidently, $U_r(x)$ is a polynomial of degree $r$ in
$x$ with integer coefficients. For example, $U_0(x)=1$,
$U_1(x)=2x$, $U_2(x)=4x^2-1$, and in general,
$U_r(x)=2xU_{r-1}(x)-U_{r-2}(x)$. Chebyshev polynomials were
invented for the needs of approximation theory, but are also
widely used in various other branches of mathematics, including
algebra, combinatorics, and number theory (see \cite{Ri}).

For $k\gs1$ we define $R_k(x)$ by
$$
R_k(x)=\frac{2tU_{k-1}(t)}{U_{k}(t)}, \qquad t=\frac1{2\sqrt{x}}.
$$
For example, $R_1(x)=1$, $R_2(x)=\frac{1}{1-x}$, and
$R_3(x)=\frac{1-x}{1-2x}$. It is easy to see that for any $k$,
$R_k(x)$ is rational in $x$. These rational functions arise in
some of the results below.

\subsection*{\bf 1.3. Preliminaries}
Apparently, for the first time the relation between restricted
permutations and Chebyshev polynomials was discovered  by Chow and
West in \cite{CW}. The main result of \cite{CW} can be formulated
as follows.

\begin{theorem}\label{tcw} {\rm (\cite[Th.~3.1]{CW})} Let $T_1=\{321,
(2,3,\dots,k,1)\}$, $T_2=\{132, (1,2,\dots,$ $k)\}$, and
$T_3=\{132, (2,3,\dots,k,1)\}$, then {\rm:}

{\rm (i)} $F_{T_1}(x)=R_k(x);$

{\rm (ii)} $F_{T_2}(x)=R_k(x);$

{\rm (iii)} $F_{T_3}(x)=R_k(x)$.
\end{theorem}

The original proof is based on the use of transfer matrices (see
Sec.~2 below). Several different proofs of various parts of this
theorem appeared recently in \cite[Th.~9]{Kr} (part (i)),
\cite[Th.~3.1]{MV1} and \cite[Th.~2]{Kr} (part (ii)), and
\cite[Th.~6]{Kr} (part (iii)). They all are based on the use of
continued fractions; the latter approach to restricted
permutations was initiated in \cite{RWZ} and developed in
\cite{MV1, JR, Kr}. In fact, there are two different ways to use
continued fractions. One of them is more geometrical, and is based
on the relation between continued fractions and Dyck paths
discovered by Flajolet (see Sec.~3 below). Another is more
analytical, and is based on the study of block decompositions (see
Sec.~4 below).

In this paper we describe these three approaches, and present
several ways to extend and generalize the above result. All our
results deal with multiple (mainly, double) restrictions, of which
one belongs to $\SS_3$ and others to $\SS_k$. Some of the results
have been published previously, and we state them here without a
proof. Observe that modulo standard symmetry operations
(complement, reversal, and inversion), there are only two
nonequivalent patterns in $\SS_3$; we choose them to be $132$ and
$321$.

For the sake of brevity, we denote by $[k]$ the identity pattern
$(1,2,\dots,k)\in\SS_k$, and by $[k,m]$ the {\it two-layered\/}
pattern $(m+1,m+2,\dots,k,1,2,\dots,m) \in\SS_k$. In general,  we
say that $\tau\in \SS_k$ is a {\it layered\/} pattern if it can be
represented as $\tau=(\tau^0,\tau^1,\dots,\tau^r)$, where each of
$\tau^i$ is a nonempty permutation of the form
$\tau^i=(m_{i+1}+1,m_{i+1}+2,\dots,m_i)$ with
$k=m_0>m_1>\dots>m_r>m_{r+1}=0$; in this case we denote $\tau$ by
$[m_0,\dots,m_r]$. Observe that our definition slightly differs
from the one used in \cite{Bo, MV3}: their layered patterns are
exactly the complements of our layered patterns.

The authors are grateful to C.~Krattenthaler and to an anonymous referee
for valuable comments that helped us to improve the presentation.
%=========================================================================
\setcounter{section}{2} \setcounter{subsection}{0}
\setcounter{theorem}{0}
\section*{\bf 2. Transfer matrices}

The main idea behind the transfer matrix approach can be described
as follows (see \cite[Th.~4.7.2]{St}). Consider a directed
multigraph on $n$ vertices $v_1,\dots,v_n$, and let $A$ denote its
weighted adjacency matrix, that is, $a_{ij}$ is the number of
edges directed from $v_i$ to $v_j$. Then the generating function
for the number of walks from $v_r$ to $v_s$ is given by
$$\frac{(-1)^{r+s}\det(I-xA;r,s)}{\det(I-xA)}\eqno(2.1)$$
where $I$ is the identity matrix and $\det(B;r,s)$ is the minor of
$B$ with the $r$th row and $s$th column deleted.

To apply this approach, one has to construct a bijection between
the permutations in question and walks in an appropriate directed
graph. We describe below two bijections of this type: the first
based on generating trees, and the second based on Dyck paths.

\subsection*{\bf 2.1. Generating trees}
Following \cite{W}, a {\it generating tree\/} is a
rooted labeled tree with the property that if $v_1$ and $v_2$ are
any two nodes with the same label and $l$ is any label, then $v_1$
and $v_2$ have exactly the same number of children with the label
$l$. To specify a generating tree it therefore suffices to
specify:
\begin{enumerate}
\item   the label of the root, and

\item   a set of succession rules explaining how to derive from the
    label of a parent the labels of all of its children.
\end{enumerate}

\begin{example}\label{aex1} {\rm(The complete binary tree)}
\rm Since all the nodes in the complete binary tree are similar,
it is enough to use only one label, which we choose to be $2$. So
we get the following description:
\begin{itemize}
\item[]  {\bf Root}: $(2)$

\item[]  {\bf Rule}: $(2)\rightarrow(2)(2)$.
\end{itemize}
\end{example}

\begin{example}\label{aex2} {\rm (The Fibonacci tree)}
\rm Here we have nodes of two different types, so we use two
labels: $1$ for a non-breeding pair and $2$ for a breeding pair.
We thus get:
\begin{itemize}
\item[] {\bf Root}: $(1)$

\item[] {\bf Rules}: $(1)\rightarrow(2)$,\quad $(2)\rightarrow(1)(2)$
\end{itemize}
\end{example}

Given a generating tree, one assigns to it a directed graph whose
vertices correspond to labels and edges from $l_i$ to $l_j$
correspond to the occurrences of $l_j$ in the succession rule
$(l_i)\rightarrow\cdots$. The graphs corresponding to the above
two examples are shown in Figure $1$.

\begin{center}
\begin{figure}[h]
\epsfxsize=2.5in
\epsffile{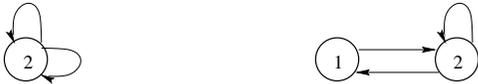} \caption{Directed graphs
for the complete binary tree and the Fibonacci tree}
\end{figure}
\end{center}

Given a permutation $\tau$, one defines a rooted tree as follows.
The nodes on level $n$ are precisely the elements of
$\SS_n(\tau)$. The parent of a permutation
$\pi=(\pi_1,\pi_2,\dots,\pi_n)$ is the unique permutation
$\pi'=(\pi_1,\dots,\pi_{j-1},\pi_{j+1},\dots,\pi_n)$ such that
$\pi_j=n$. We denote the resulting tree $\T(\tau)$. Similarly, the
tree corresponding to the set $\SS_n(T)$ is denoted by $\T(T)$.

Chow and West \cite{CW} proved that the succession rules for the
tree $\T(123,(k-1,\dots,1,k))$ are
$$\begin{array}{ll}
(l)\rightarrow(2)\cdots(l)(l+1),   & \quad l< k-1\\[4pt]
(k-1)\rightarrow(2)\cdots(k-1)(k-1),&
\end{array}$$
and the label of the root is $(2)$. The corresponding graph is
shown in Figure $2$.

\begin{center}
\begin{figure}[h]
\epsfxsize=2.5in \epsffile{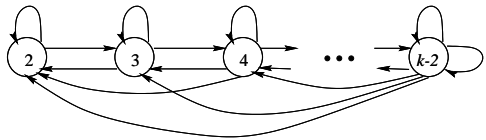} \caption{The directed graph
for $\T(123,(k-1,\dots,2,1,k))$ }
\end{figure}
\end{center}

The corresponding transfer matrix is
$$A_k=\left[
    \begin{array}{llllllll}
    1   &   1   &   0   &0 & \cdots    &   0   &   0   &   0\\
    1   &   1   &   1   &0 & \cdots    &   0   &   0   &   0\\
    1   &   1   &   1   &1 & \cdots    &   0   &   0   &   0\\
    \vdots  &   \vdots  &\vdots& \vdots&\ddots &   \vdots  &   \vdots  &   \vdots\\
    1   &   1   &   1   &1 & \cdots    &   1   &   0   &   0\\
    1   &   1   &   1   &1 & \cdots    &   1   &   1   &   0\\
    1   &   1   &   1   &1 & \cdots    &   1   &   1   &   1\\
    1   &   1   &   1   &1 & \cdots    &   1   &   1   &   2\\
    \end{array}
    \right]$$
Besides, Chow and West proved that the graphs, and hence the
transfer matrices for $\T(213,(k,1,2\dots,k-1))$ and
$\T(213,(1,2,\dots,k))$ are exactly the same (though the
succession rules may vary). The number of permutations in
$\SS_n(T)$ (in all cases) is thus equal to the number of walks of
length $n$ starting from the vertex $2$. However, since each
vertex is connected to vertex $2$ by exactly one edge, this number
is equal to the number of walks of length $n+1$ starting at vertex
$2$ and ending at the same vertex. The generating function for
this number is given by $(2.1)$ with $A=A_k$. It is proved in
\cite{CW} that the determinants in question satisfy linear
recurrences of order two very similar to that for Chebyshev
polynomials, which almost immediately yields Theorem \ref{tcw},
since $T_1=c(\{123, (k-1,\dots,2,1,k)\})$,
$T_2=r\circ c(\{213, (1,2,\dots,k)\})$, and
$T_3=r\circ c(\{213, (k,1,2,\dots,,k-1)\})$.

\subsection*{\bf 2.2. Dyck paths}
A {\em Dyck path\/} is a path in the plane integer lattice ${\bf
Z}^2$, consisting of up-steps $(1,1)$ and down-steps $(1,-1)$,
which never passes below the $x$-axis.

Following \cite{Kr}, we define a bijection $\Phi$ between
permutations in $\SS_n(132)$ and Dyck paths from the origin to the
point $(2n,0)$. Let $\pi=(\pi_1,\dots,\pi_n)$ be a $132$-avoiding
permutation. We read the permutation $\pi$ from left to right and
successively generate a Dyck path. When $\pi_j$ is read, then in
the path we adjoint as many up-steps as necessary, followed by a
down-step from height $h_j+1$ to height $h_j$ (measured from the
$x$-axis), where $h_j$ is the number of elements in
$\pi_{j+1},\pi_{j+2},\dots,\pi_n$ which are larger that $\pi_j$.

For example, let $\pi=534261$. The first element to be read is
$5$. There is one element in $34261$ which is larger than $5$,
therefore the path starts with two up-steps followed by a
down-step, thus reaching height $1$. Next $3$ is read. There are
$2$ elements in $4261$ which are larger than $3$, therefore the
path continues with two up-steps followed by a down-step, thus
reaching height $2$. Etc.

Conversely, given a Dyck path starting at the origin and returning
to the $x$-axis, the obvious inverse of the bijection $\Phi$
produces a $132$-avoiding permutation.

It is proved in \cite{Kr} that the bijection $\Phi$ takes
permutation in $\SS_n(132,[k])$ to Dyck paths that never pass
above the line $y=k-1$. Evidently, such paths correspond
bijectively to walks of length $2n$ starting at vertex $0$ in the
graph shown in Figure $3$.

\begin{center}
\begin{figure}[h]
\epsfxsize=2.5in \epsffile{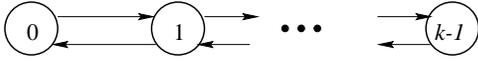} \caption{The directed graph
for Dyck paths in a strip}
\end{figure}
\end{center}

Using again $(2.1)$ one gets Theorem \ref{tcw}(ii). A further
study of the bijection $\Phi$ yields part $(iii)$ of the same
Theorem.
%======================================================================
%======================================================================
\setcounter{section}{3} \setcounter{subsection}{0}
\setcounter{theorem}{0}
\section*{\bf 3. Continued fractions}
The relation between restricted permutations and continued
fractions was discovered by Robertson, Wilf, and Zeilberger in
\cite{RWZ}. The main result in \cite{RWZ} can be formulated as
follows. Let $G_\tau^r(x)$ be the generating function for the
number of permutations in $\SS_n(132)$ containing a pattern $\tau$
exactly $r$ times.

\begin{theorem}{\rm (Robertson, Wilf, and Zeilberger \cite[Th.~1]{RWZ})}\label{thrwz}
$$\sum_{r\gs 0} G_{123}^r(x)z^r=\frac{1}
        {1-\dfrac{xz^{\binom{0}{2}}}
        {1-\dfrac{xz^{\binom{1}{2}}}
        {1-\dfrac{xz^{\binom{2}{2}}}
        {\ddots}}}}$$
in which the $j$th numerator is $xz^{\binom{j-1}{2}}$.
\end{theorem}

To prove this, let $\pi$ be a permutation avoiding $132$. Then
each letter in $\pi$ to the left of $n$ must be greater than any
letter to the right of $n$. Thus, if $\pi=(\pi',n,\pi'')$ (where
both $\pi'$ and $\pi''$ must necessarily be $132$-avoiding), then
    $$(123)\pi=(123)\pi'+(12)\pi'+(123)\pi'',$$
where $(\tau)\pi$ is the number of occurrences of $\tau$ in $\pi$.
It follows that the generating function
    $$F(x,y,z)=\sum_{\pi\in\SS(132)} x^{(1)\pi}y^{(12)\pi}z^{(123)\pi}$$
satisfies the equation $F(x,y,z)=1+xF(xy,yz,z)F(x,y,z)$.
Equivalently,
    $$F(x,y,z)=\frac{1}{1-xF(xy,yz,z)},$$
and the theorem follows by induction after plugging in $y=1$.\\

This result was generalized by Mansour and Vainshtein \cite{MV1},
by Krattenthaler \cite{Kr}, and by Jani and Rieper \cite{JR} to
the case of permutations containing the pattern $[k]=12\dots k$
exactly $r$ times. It turns out that
$$\sum_{r\geq 0} G_{[k]}^r(x)z^r
=\frac1{1-\dfrac{xz^{d_1}}{1-\dfrac{xz^{d_2}}{1-\dfrac{xz^{d_3}}{\dots}}}},
\eqno(3.1)
$$
where $d_j=\binom{j-1}{k-1}$.\\

The proof in \cite{MV1} is a straightforward generalization of the
above proof of Theorem \ref{thrwz}.  The proof in \cite{Kr} is
based on the bijection $\Phi$ between $132$-avoiding permutations
and Dyck paths described in the previous section and on the result
of Flajolet \cite[Th.~1]{F} presenting the generating function for
the Dyck paths in terms of continued fractions. The proof of
\cite{JR} is based on a bijection between $132$-avoiding
permutations and rooted ordered trees, which can be obtained from the
bijection $\Phi$ via the standard bijection between rooted ordered
trees and Dyck paths through a depth-first traversal of the trees
(see \cite[Prop.~6.2.1, Cor.~6.2.3]{St2}). For further generalizations
and interesting combinatorial applications see \cite{BCS}.

It was observed in \cite{MV1} that $R_k(x)$ is the $k$th
approximant for the continued fraction
    $$\dfrac{1}{1-\dfrac{x}{1-\dfrac{x}{1-\ddots_{{\ }_{\ } }}}}$$
so $(3.1)$ for $r=0$ immediately gives Theorem \ref{tcw}(ii).

Paper \cite{Kr} contains the description of a bijection $\Psi$
between $123$-avoiding permutations and Dyck paths. This
bijection, combined with Roblet and Viennot's continued fraction
representation of the generating function for Dyck paths
\cite[Prop.~1]{RV} gives the first part of Theorem \ref{tcw}.
%=========================================================================
\setcounter{section}{4} \setcounter{subsection}{0}
\setcounter{theorem}{0}
\section*{\bf 4. Block decompositions}

The core of this approach initiated by Mansour and Vainshtein
\cite{MV2} lies in the study of the structure of $132$-avoiding
permutations, and permutations containing a given number of
occurrences of $132$.

Let us start with the simplest case of $132$-avoiding
permutations. It was noticed in \cite{MV2} that if
$\alpha\in\SS_n(132)$ and $\alpha_t=n$, then
$\alpha=(\alpha',n,\alpha'')$ where $\alpha'$ is a permutation of
the numbers $n-t+1, n-t+2,\dots,n$, $\alpha''$ is a permutation of
the numbers $1,2,\dots,n-t$, and both $\alpha'$ and $\alpha''$
avoid $132$. This representation is called the block decomposition
of $\alpha$, see Figure $4$.

\begin{center}
\begin{figure}[h]
\epsfxsize=2.5in \epsffile{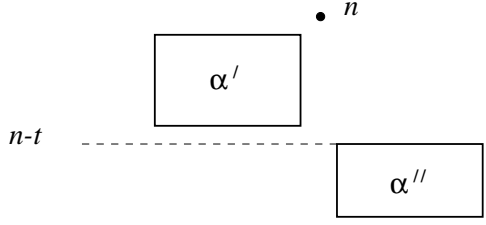} \caption{The block
decomposition for $\alpha\in\SS_n(132)$}
\end{figure}
\end{center}

The simple observation allows to formulate a general result
concerning permutations avoiding $132$ and arbitrary pattern
$\tau=(\tau_1,\dots,\tau_k)\in\SS_k(132)$. Recall that $\tau_i$ is
said to be a {\it right-to-left maximum\/} if $\tau_i>\tau_j$ for
any $j>i$. Let $m_0=k,m_1,\dots,m_r$ be the right-to-left maxima
of $\tau$ written from left to right. Then $\tau$ can be
represented as
    $$\tau=(\tau^0,m_0,\tau^1,m_1,\dots,\tau^r,m_r),$$
where each of $\tau^i$ may be possibly empty, and all the
entries of $\tau^i$ are greater than all the entries of
$\tau^{i+1}$. Define the $i$th {\it prefix\/} of $\tau$ by
$\pi^i=(\tau^0,m_0,\dots,\tau^i,m_i)$ for $1\ls i\ls r$ and
$\pi^0=\tau^0$, $\pi^{-1}=\varnothing$. Besides, the $i$th {\it
suffix\/} of $\tau$ is defined by
$\sigma^i=(\tau^i,m_i,\dots,\tau^r,m_r)$ for $0\ls i\ls r$ and
$\sigma^{r+1}=\varnothing$.

\begin{theorem}\label{a1} {\rm (\cite[Th.~1]{MV2})}
For any $\tau\in\SS_k(132)$, $F_\tau(x)$ is a rational function
satisfying the relation
$$F_\tau(x)=1+x\sum_{j=0}^r\bigl(F_{\pi^j}(x)-F_{\pi^{j-1}}(x)\bigr)
F_{\sigma^j}(x).$$
\end{theorem}

The proof is rather straightforward. Let
$\alpha=(\alpha',n,\alpha'')$ be the block decomposition of
$\alpha\in\SS_n(132)$. It is easy to see that $\alpha$ contains
$\tau$ if and only if there exists $i$, $0\ls i\ls r+1$, such that
$\alpha'$ contains $\pi^{i-1}$ and $\alpha''$ contains $\sigma^i$.
Therefore, $\alpha$ avoids $\tau$ if and only if there exists $i$,
$0\ls i\ls r$, such that $\alpha'$ avoids $\pi^i$ and contains
$\pi^{i-1}$, while $\alpha''$ avoids $\sigma^i$. We thus get the
following relation:
   $$f_\tau(n)=\sum_{t=1}^n\sum_{j=0}^r f_{\pi^j}^{\pi^{j-1}}(t-1)f_{\sigma^j}(n-t),$$
where $f_\tau(n)=|\SS_n(132,\tau)|$, and $f_{\tau}^{\rho}(n)$ is
the number of permutations in $\SS_n(\tau)$ containing $\rho$ at
least once. To obtain the recursion for $F_\tau(x)$ it remains to
observe that
    $$f_{\pi^j}^{\pi^{j-1}}(l)+f_{\pi^{j-1}}(l)=f_{\pi^j}(l)$$
for any $l$ and $j$, and to pass to generating functions.
Rationality of $F_\tau(x)$ follows easily by induction.\\

Theorem \ref{a1} allows to reduce the calculation of $F_\tau(x)$
to finding similar functions for several simpler patterns. For
example, if $\tau=[k,m]$ is a two-layered pattern, then Theorem
\ref{a1} gives
$$F_{[k,m]}(x)=1+xF_{[k-m-1]}(x)F_{[k,m]}(x)+x(F_{[k,m]}(x)-F_{[k-m-1]}(x))F_{[m]}(x)\eqno(4.1)$$
and hence $F_{[k,m]}(x)$ can be expressed via $F_{[m]}(x)$ and
$F_{[k-m-1]}(x)$, (see Theorem \ref{a51} below).\\

Consider now the case of permutations containing $132$ exactly
once.

\begin{center}
\begin{figure}[h]
\epsfxsize=2.5in \epsffile{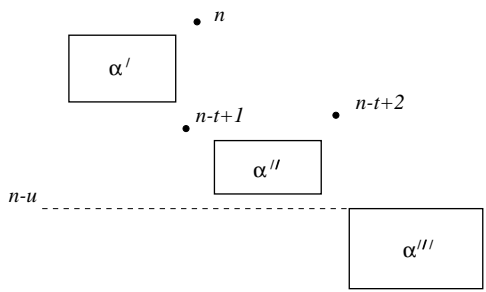} \caption{The block
decomposition for $\alpha\in\SS_n$ containing $132$ exactly once}
\end{figure}
\end{center}

\begin{theorem}\label{thbd1}
Let $\alpha\in\SS_n$ contain $132$ exactly once. Then the block
decomposition of $\alpha$ can have one of the following three
forms:

(i) there exists $t$ such that $\alpha=(\alpha',n,\alpha'')$,
where $\alpha'$ is a permutation of $n-1,n-2,\dots,n-t+1$
containing $132$ exactly once, and $\alpha''$ is a permutation of
$1,2,\dots,n-t$ avoiding $132$;

(ii) there exists $t$ such that $\alpha=(\alpha',n,\alpha'')$,
where $\alpha'$ is a permutation of $n-1,n-2,\dots,n-t+1$ avoiding
$132$, and $\alpha''$ is a permutation of $1,2,\dots,n-t$
containing $132$ exactly once;

(iii) there exist $t,u$ such that
$\alpha=(\alpha',n-t+1,n,\alpha'',n-t+2,\alpha''')$, where
$\alpha'$ is a permutation of $n-1,n-2,\dots,n-t+3$ avoiding
$132$, $\alpha''$ is a permutation of $n-t,n-t-1,\dots,n-u+1$
avoiding $132$, and $\alpha'''$ is a permutation of
$1,2,\dots,n-u$ avoiding $132$.
\end{theorem}
\begin{proof} Let $\alpha\in\SS_n$ contain $132$ exactly once.
There are two possibilities: either the only occurrence of $132$
in $\alpha$ does not contain $n$, or it contains $n$. In the first
case we get immediately that any entry of $\alpha$ to the right of
$n$ is less than any entry of $\alpha$ to the left of $n$, since
otherwise one gets an occurrence of $132$ involving $n$. In the
second case, let $i\; n\; j$ be the occurrence of $132$ in
$\alpha$. First, we have $j=i+1$, since otherwise either $(i+1)\; n
\;j$ or $i\; n\; (i+1)$ would be a second occurrence of $132$ in
$\alpha$. Next, $i$ immediately precedes $n$ in $\alpha$, since if
$l$ lies between $i$ and $n$, then either $l\; n\; j$ or $i\;l\;j$
would be a second occurrence of $132$ in $\alpha$. Finally, any
entry to the right of $j$ is less than any entry between $n$ and
$j$, which in turn, is less than any entry to the left of $i$
(the proof is similar to the analysis of the first case).
\end{proof}

The block decompositions of types $(i)$ and $(ii)$ are similar to
that for $\alpha\in\SS_n(132)$ (see Figure 4). The block decomposition of type
$(iii)$ is shown in Figure $5$.

For a far reaching generalization of this idea, allowing to
enumerate permutations with any given number of occurrences of
$132$, see \cite{MV4}.
%=========================================================================
\section*{\bf 5. Counting $132$-restricted permutations}
\setcounter{section}{5} \setcounter{subsection}{0}
\setcounter{theorem}{0}

Throughout this section we write $\bar{F}_\tau(x)$ instead of
$F_{132,\tau}(x)$.

\subsection*{\bf 5.1. Avoiding $132$ and another
pattern}\addtocounter{subsection}{1} By Theorem \ref{tcw}~(ii)
and~(iii), the pairs $\{132,[k]\}$ and $\{132,[k,1]\}$ are {\it
Wilf-equivalent}, that is $\bar F_{[k]}(x)=\bar F_{[k,1]}(x)$. It
turns out that this equivalence class can be extended as follows.

\begin{theorem} {\rm (\cite[Th.~2.4]{MV2})}\label{a51} For any $m$,
$1\leq m\leq k-1$,
$$
\bar F_{[k,m]}(x)=R_k(x).
$$
\end{theorem}

This result follows immediately from Theorem \ref{tcw} and
$(4.1)$. \\

A further extension is provided by the following result. We say
that $\tau\in \SS_k$ is a {\it wedge\/} pattern if it can be
represented as $\tau=(\tau^1,\rho^1,\dots,\tau^r,\rho^r)$ so that
each of $\tau^i$ is nonempty, $(\rho^1,\rho^2,\dots,\rho^r)$ is a
layered permutation of $1,\dots,s$ for some $s$, and
$(\tau^1,\tau^2,\dots,\tau^r)=(s+1,s+2,\dots,k)$. For example,
$645783912$ is a wedge pattern; here $r=3$, $s=5$, $\tau^1=(6)$,
$\tau^2=(7,8)$, $\tau^3=(9)$, $\rho^1=(4,5)$, $\rho^2=(3)$,
$\rho^3=(1,2)$. Evidently, $[k,m]$ is a wedge pattern for any $m$.

\begin{theorem} {\rm(\cite[Th.~2.6]{MV2})}\label{a52} $\bar
F_{\tau}(x)=R_k(x)$ for any wedge pattern $\tau\in \SS_k(132)$.
\end{theorem}

This result follows easily from Theorem \ref{a1}. The proof of the
above two results are purely analytical. It would be very
interesting to find a bijective proof of these results.

The case of general layered patterns can be also expressed in
terms of Chebyshev polynomials, using the same technique of block
decompositions. Since the expressions become rather cumbersome, we
present here only the case of a $3$-layered pattern.

\begin{theorem} {\rm (\cite[Th.~2.5]{MV2})}\label{a53}
 For any $k>m_1>m_2>0$,
$$
\bar F_{[k,m_1,m_2]}(x)=2t\frac{U_{a+b}(t) U_{a+c-1}(t)U_{b+c}(t)+
U_{b-1}(t)U_b(t)} {U_{a+b}(t) U_{a+c}(t)U_{b+c}(t)},\qquad
t=\frac1{2\sqrt{x}},
$$
where $a=k-m_1$, $b=m_1-m_2$, $c=m_2$.
\end{theorem}

\subsection*{\bf 5.2. Avoiding $132$ and containing another pattern}\addtocounter{subsection}{1}
Denote by $G^r_{\tau}(x)$ the generating function for the number
of permutations in $\SS_n(132)$ containing a pattern $\tau\in
\SS_k$ exactly $r$ times. The continued fraction representation
$(3.1)$ for $r=1$ immediately gives the following result.

\begin{theorem} {\rm (\cite[Th.~3.1]{MV1})}\label{b54} For any
$k\geq1$,
$$
G^1_{[k]}(x)=\frac1{U^{2}_{k}(t)},\qquad t=\frac1{2\sqrt{x}}.
$$
\end{theorem}

This result may be extended in two directions. First, let us fix
$r=1$ and consider other patterns $\tau$. The case of $\tau=[k,1]$
was studied in \cite{Kr}.

\begin{theorem} {\rm (\cite[Th.~7]{Kr})}\label{b55} For $k\geq2$,
$$
G^1_{[k,1]}(x)=\frac{1}{4t^2U_{k-2}(t)U_{k}(t)}, \qquad
t=\frac1{2\sqrt{x}}.
$$
\end{theorem}

The proof involves  bijection $\Phi$ and the result of Flajolet
mentioned in Section $3$.

A more general case of $\tau=[k,m]$ is investigated in \cite{MV2},
based on block decompositions.

\begin{theorem} {\rm (\cite[Th.~3.4]{MV2})}\label{b56} For any
$k>m>0$,
$$
G^1_{[k,m]}(x)=\frac{1}{2tU_{k}(t)U_{m}(t)U_{k-m-1}(t)}, \qquad
t=\frac1{2\sqrt{x}}.
$$
\end{theorem}

Another direction would be to increase the value of $r$. A
generalization of Theorem \ref{b54} to the case $1\leq r\leq k$
was obtained in \cite{MV1}, directly from $(3.1)$.

\begin{theorem} {\rm (\cite[Th.~3.1]{MV1})}\label{b57} For any $r$,
$1\leq r \leq k$,
$$
G^r_{[k]}(x)=\frac{U^{r-1}_{k-1}(t)} {(2t)^{r-1}U^{r+1}_{k}(t)},
\qquad t=\frac1{2\sqrt{x}}.
$$
\end{theorem}

A similar generalization of Theorem \ref{b55} is as follows.

\begin{theorem} {\rm (\cite[Th.~7]{Kr})}\label{B58} For any $r$,
$1\leq r \leq k-1$,
$$
G^r_{[k,1]}(x)=
%\\
\frac1{U_{k-3}(t)U_{k}(t)}\sum_{l|r}\frac1{l+1}\binom{2l}{l}
(2t)^{1-2l-\frac{r}{l}}\left(\frac{U_{k-3}(t)}
{U_{k-2}(t)}\right)^{r/l}, \qquad t=\frac1{2\sqrt{x}}.
$$
\end{theorem}
The ideas behind the proof are the same as in the proof of Theorem
\ref{b55}.

Theorem \ref{b57} can be extended to cover a wider range of $r$'s.

\begin{theorem} {\rm (\cite[Th.~4.1]{MV1})}\label{b59} For any $r$,
$1\leq r \leq k(k+3)/2$,
$$
G^r_{[k]}(x)= \frac{U^{r-1}_{k-1}(t)}{(2t)^{r-1}U^{r+1}_{k}(t)}
\sum_{j=0}^{\lfloor(r-1)/k\rfloor}
\binom{r-kj+j-1}{j}\left(\frac{(2t)^{\frac{k-2}{k}}U_{k}(t)}
{U_{k-1}(t)}\right)^{kj},
$$
where $t=1/2\sqrt{x}$.
\end{theorem}

The case of general $r$ was treated in \cite{Kr}.

\begin{theorem} {\rm (\cite[Th.~3]{Kr})}\label{b510} For any $r$,
$$
G^r_{[k]}(x)=
\sum\left(\binom{l_1+l_2-1}{l_2}\binom{l_2+l_3-1}{l_3}
\dots\right)\frac{U_{k-1}^{l_1-1}(t)} {U_{k}^{l_1+1}(t)}
(2t)^{-(l_1-1)-2(l_2+l_3+\dots)},
$$
where the sum is over all nonnegative integers $l_1,l_2,\dots$
with
$$
l_1\binom{k-1}{k-1}+l_2\binom{k}{k-1}+l_3\binom{k+1}{k-1}+\dots=r,
$$
and $t=1/2\sqrt{x}$.
\end{theorem}

\subsection*{\bf 5.3. Containing $132$ exactly once and avoiding another
pattern}\addtocounter{subsection}{1} Denote by $H_{\tau}(x)$ the
generating function for the number of permutations in $\SS_n$
avoiding a pattern $\tau\in \SS_k$ and containing $132$ exactly
once. We start from the following result obtained in \cite{MV1}.

\begin{theorem} {\rm (\cite[Th.~4.2]{MV1})}\label{c511} For any
$k\geq3$,
$$
H_{[k]}(x)=\frac{1}{4t^2U_{k}^2(t)} \sum_{j=1}^{k-2}U_{j}^2(t),
\qquad t=\frac1{2\sqrt{x}}.
$$
\end{theorem}

The idea behind the proof is similar to that of the proof of
Theorem \ref{thrwz} explained above in Section $3$.

This result can be extended to the case of general $2$-layered
patterns as follows.

\begin{theorem} {\rm (i)} For any $k\geq4$,
$$
H_{[k,1]}(x)=\frac{1}{4t^2U_{k}^2(t)}
\left(\sum_{j=1}^{k-2}U_{j}^2(t)-1\right), \qquad
t=\frac1{2\sqrt{x}};
$$
besides,
$$
H_{[3,1]}(x)=\frac{x^3}{1-2x}.
$$

{\rm (ii)} For any $k\geq5$,
$$
H_{[k,2]}(x)=\frac{1}{4t^2U_{k}^2(t)}
\left(\sum_{j=1}^{k-2}U_{j}^2(t)-\frac{2tU_{k-3}(t)}{U_{k-2}(t)}-2\right),
\qquad t=\frac1{2\sqrt{x}};
$$
besides,
$$
H_{[4,2]}(x)=\frac{x^3(1+x)}{(1-x)(1-3x+x^2)}.
$$

{\rm (iii)} For any $k\geq6$ and any $m$, $3\leq m\leq k/2$,
$$\begin{array}{l}
H_{[k,m]}(x)=\frac{1}{4t^2U_{k}^2(t)}
\bigg(\sum_{j=1}^{k-m-2}U_{j}^2(t)+
\sum_{j=1}^{m-1}U_{j}^2(t)-1\\
\quad\quad\quad\quad\quad\quad\quad\quad\quad\quad\quad+U_{k-1}(t)U_{m-1}(t)U_{k-m-2}(t)\bigg),
\qquad t=\frac1{2\sqrt{x}}.
\end{array}
$$
\end{theorem}
\noindent{\bf Remark}. Clearly, the cases $k=3$, $m=2$; $k=4,5$,
$m=3$; $k\geq6$, $m>k/2$ are also covered, since $H_{[k,m]}(x)=
H_{[k,k-m]}(x)$.
\begin{proof}
By Theorem \ref{thbd1}, we have exactly three possibilities for
the block decomposition of an arbitrary $\alpha\in\SS_n$. Let us
write an equation for $H_{[k,1]}(x)$ with $k\geq4$. The
contribution of the first decomposition above is
$xH_{[k-2]}(x)\bar
F_{[k,1]}(x)+x\big(H_{[k,1]}(x)-H_{[k-2]}(x)\big)$. Here the first
term corresponds to the case $\alpha'$ avoids $[k-2]$ and
$\alpha''$ avoids $[k,1]$, while the second term corresponds to
the case $\alpha'$ avoids $[k,1]$ but contains $[k-2]$, and
$\alpha''=\varnothing$.

The contribution of the second possible decomposition is $x\bar
F_{[k-2]}(x)H_{[k,1]}(x)$; here $\alpha''$ contains $132$, and
hence is always distinct from $\varnothing$.

Finally, the contribution of the third possible decomposition is
$$x^3\bar F_{[k-2]}^2(x)\bar F_{[k,1]}(x)+ x^3\bar
F_{[k-2]}(x)(\bar F_{[k,1]}(x)-\bar F_{[k-2]}(x)).$$

Here the first
term corresponds to the case $\alpha',\alpha''$ avoid $[k-2]$,
$\alpha'''$ avoids $[k,1]$, while the second term corresponds to
the case $\alpha'$ avoids  $[k-2]$, $\alpha''$ avoids $[k,1]$  but
contains $[k-2]$, and $\alpha'''=\varnothing$.

Solving the obtained linear equation and using Theorems
\ref{tcw},~\ref{c511}, and well known identities involving
Chebyshev polynomials (see e.g. \cite[Lem.~4.1]{MV2}), we get the
desired expression for $H_{[k,1]}(x)$, $k\geq4$.

In the case $k=3$, the contributions of the first and the second
decompositions degenerate to $xH_{[3,1]}(x)$ each, while the
contribution of the third decomposition degenerates to $x^3$
(which means that the only permutation having this decomposition
is $132$ itself). The result follows immediately.

Let us consider now the case of $H_{[k,2]}(x)$ with $k\geq5$. The
contribution of the first decomposition is $xH_{[k-3]}(x)\bar
F_{[k,2]}(x)+x\big(H_{[k,2]}(x)- H_{[k-3]}(x)\big)\bar
F_{[2]}(x)$. Here the first term corresponds to the case $\alpha'$
avoids $[k-3]$ and $\alpha''$ avoids $[k,2]$, while the second
term corresponds to the case $\alpha'$ avoids $[k,2]$ but contains
$[k-3]$, and $\alpha''$ avoids $[2]$.

The contribution of the second decomposition is $x\bar
F_{[k-3]}(x)H_{[k,2]}(x)$.

Finally, the contribution of the third decomposition is
$$\begin{array}{l}
x^3\bar F_{[k-3]}^2(x)\bar F_{[k,2]}(x)+ x^3\bar
F_{[k-3]}(x)\big(\bar F_{[k,2]}(x)-\bar F_{[k-3]}(x)\big) \bar
F_{[2]}(x)\\
\quad\quad\quad\quad\quad\quad\quad\quad\quad\quad\quad\quad\quad+x^3\big(\bar
F_{[k-2]}(x)-\bar F_{[k-3]}(x)\big) \bar F_{[2]}(x).\end{array}
$$
Here the first term corresponds to the case $\alpha',\alpha''$
avoid $[k-3]$, $\alpha'''$ avoids $[k,2]$, the second term
corresponds to the case $\alpha'$ avoids  $[k-3]$, $\alpha''$
avoids $[k,2]$  but contains $[k-3]$, and $\alpha'''$ avoids
$[2]$, while the third term corresponds to the case $\alpha'$
avoids $[k-2]$ but contains $[k-3]$, $\alpha''=\varnothing$,
 and $\alpha'''$ avoids $[2]$.

The expression for $H_{[k,2]}(x)$ follows easily from this
equation and Theorems \ref{tcw} and \ref{c511}.

The case $k=4$ is treated similarly.

For general $m$, $3\leq m\leq k-m$, the contribution of the first
decomposition equals $xH_{[k-m-1]}(x)\bar
F_{[k,m]}(x)+x\big(H_{[k,m]}(x)-H_{[k-m-1]}(x)\big) \bar
F_{[m]}(x)$, the contribution of the second decomposition equals
$x\bar F_{[k-m-1]}(x)H_{[k,m]}(x)+ x\big(\bar F_{[k,m]}(x)-\bar
F_{[k-m-1]}(x)\big)H_{[m]}(x)$, and the contribution of the third
structure decomposition
$$\begin{array}{l}
x^3\bar F_{[k-m-1]}^2(x)\bar F_{[k,m]}(x)+x^3\bar F_{[k-m-1]}(x)
\big(\bar F_{[k,m]}(x)-\bar F_{[k-m-1]}(x)\big)\bar F_{[m]}(x)\\
\quad\quad\quad\quad\quad\quad\quad\quad\quad\quad\quad\quad\quad+x^3\big(\bar
F_{[k,m]}(x)-\bar F_{[k-m-1]}(x)\big)\bar F_{[m-1]}(x) \bar
F_{[m]}(x). \end{array}
$$
The final result again follows from Theorems \ref{tcw} and
\ref{c511}.
\end{proof}

\subsection*{\bf 5.4. Containing $132$ and another pattern exactly once}\addtocounter{subsection}{1}
Denote by $\Phi_{\tau}(x)$ the generating function for the number
of permutations in $\SS_n$ containing both $132$ and a pattern
$\tau\in \SS_k$ exactly once. We start from the following result.

\begin{theorem}\label{d513} For any $k\geq1$,
$$
\Phi_{[k]}(x)=\frac{1}{4t^3U_{k}^2(t)}
\sum_{i=1}^{k-2}\frac{\sum_{j=1}^{k-i}U_{j}^2(t)-1}
{U_{k-i}(t)U_{k-i+1}(t)}, \qquad t=\frac1{2\sqrt{x}}.
$$
\end{theorem}

\begin{proof} The three possible block decompositions of permutations
containing $132$ exactly once are described in Theorem
\ref{thbd1}. Let us find the recursion for $\Phi_{[k]}(x)$. It is
easy to see that the contribution of the first decomposition
equals
$$x\Phi_{[k-1]}(x)\bar F_{[k]}(x)+xH_{[k-1]}(x)G_{[k]}^1(x),$$
the contribution of the second decomposition equals
$$xG_{[k-1]}^1(x)H_{[k]}(x)+x\bar F_{[k-1]}(x)\Phi_{[k]}(x),$$
while the contribution of the third decomposition equals
$$2x^3G_{[k-1]}^1(x)\bar F_{[k-1]}(x)\bar F_{[k]}(x)+x^3\bar F_{[k-1]}^2(x)
G_{[k-1]}^1(x).$$ Solving the obtained recursion with the initial
condition $\Phi_{[2]}(x)=0$ and using Theorems \ref{tcw},
\ref{b54}, and \ref{c511}, we get the desired result.
\end{proof}

In particular, for $k=3$ we get $\Phi_{[3]}(x)=2x^5(1-2x)^{-3}$,
which means that the number of permutations in $\SS_n$ containing
both $132$ and $123$ exactly once equals $(n-3)(n-4)2^{n-5}$ (see
\cite[Th.~1.3.18]{R}).

Similarly to the previous section, this result can be extended to
the case of general $2$-layered patterns. Since the answers become
very cumbersome, we present here only the simplest case.

\begin{theorem}\label{d514} For any $k\geq4$,
$$\begin{array}{l}
\Phi_{[k,1]}(x)=\frac {1} {8t^3U_{k}^2(t)}\left(
\frac{U_k(t)}{2t^2U_{k-2}(t)}\sum_{i=1}^{k-4}
\frac{\sum_{j=1}^{k-i-2}U_{j}^2(t)-1}{U_{k-i-2}(t)U_{k-i-1}(t)}\right.\\
\quad\quad\quad\quad\quad\quad\quad\quad\quad\quad\quad+\left.\frac{1}{2tU_{k-2}^2(t)}\sum_{i=1}^{k-4}U_{i}^2(t)+
\frac{U_{k-3}(t)}{2tU_{k-1}(t)}+\frac{U_{k}(t)}{U_{k-1}(t)}\right),
\end{array}
$$
where $t=1/2\sqrt{x}$.
\end{theorem}

\subsection*{\bf 5.5. Generalizations}\addtocounter{subsection}{1} Here we present several
directions to generalize the results of the previous sections.
The first of these directions is to consider more than one additional
restriction. For example, the following result is true. Let $\bar
F_{\tau_1,\tau_2}(x)$ be the generating function for the number of
permutations in $\SS_n(132,\tau_1, \tau_2)$. Assume that
$\tau_1=[k,m]$, $k-m\gs m$, and $\tau_2=[l]$. It is easy to see
that the only interesting case is $l>k-m$, since otherwise $\bar
F_{[k,m],[l]}(x)=\bar F_{[l]}(x)$.

\begin{theorem}\label{e515} Let $l> k-m \geq m$, then
$$
\bar F_{[k,m],[l]}(x)=R_k(x)-\big(xR_{k-m}(x)R_m(x)\big)^{l+m-k}
(R_k(x)-R_{k-m}(x)).
$$
\end{theorem}
\begin{proof} Let $\alpha\in\SS_n(132)$, then either
$\alpha=\varnothing$ (which means that $n=0$), or
$\alpha=(\alpha', n, \alpha'')$, where both $\alpha'$ and
$\alpha''$ avoid $132$. According to this dichotomy, we get the
following recursion:
$$
\bar F_{[k,m],[l]}(x)=1+x\bar F_{[k-m-1]}(x)\bar F_{[k,m],[l]}(x)+
x\big(\bar F_{[k,m],[l-1]}(x)-\bar F_{[k-m-1]}(x)\big)\bar
F_{[m]}(x).
$$
Solving this recursion with the initial condition $\bar
F_{[k,m],[k-m]}(x)= \bar F_{[k-m]}(x)$ and using Theorem
\ref{tcw}, we get the desired result.
\end{proof}

Similarly, let $G_{\tau_1;\tau_2}(x)$ be the generating function
for the number of permutations in $\SS_n(132,\tau_1)$ that contain
$\tau_2$ exactly once. As before, we assume $\tau_1=[k,m]$,
$k-m\gs m$, and $\tau_2=[l]$. Once again, the case $m\gs l$ is of
no interest, since in this case $G_{[k,m];[l]}(x)=G_{[l]}(x)$.

\begin{theorem}\label{e516} {\rm (i)} Let $l> k-m\geq m$, then
$$\begin{array}{l}
G_{[k,m];[l]}(x)=\frac{1}{U_m(t)U_{k-m}(t)}\left(\frac{U_{m-1}(t)}
{U_m(t)}\right)^{l-m}\left(\frac{U_{k-m-1}(t)}{U_{k-m}(t)}\right)^{l+m-k}\\
\quad\quad\quad\quad\quad\quad\quad\quad\times\left(\sum_{j=m+2}^{k-m}\frac{U_{j-m-2}(t)}
{U_{j-2}(t)U_{j-1}(t)}\left(\frac{U_{m-1}(t)}{U_m(t)}\right)^{m+1-j}+1\right),
\end{array}
$$
where $t=1/2\sqrt{x}$.

{\rm (ii)} Let $k-m\geq l> m$, then
$$\begin{array}{l}
G_{[k,m];[l]}(x)=\\
\quad\frac{1}{U_l(t)U_m(t)}\left(\frac{U_{m-1}(t)}{U_m(t)}\right)^{l-m}
\left(\sum_{j=m+1}^{l}\frac{U_{j-m-1}(t)}
{U_{j-1}(t)U_{j}(t)}\left(\frac{U_{m-1}(t)}{U_m(t)}\right)^{m-j}+1\right),
\end{array}
$$
where $t=1/2\sqrt{x}$.
\end{theorem}

Another possible direction is to replace $132$ by some restriction
of length $4$ or more having a similar restrictive power. Define
$L_p$ as the set of all patterns in $\SS_p$ of the form
$\pi_11\pi_22\pi_3$, where $\pi_2$ is nonempty. Evidently,
$L_3=\{132\}$; for $p=4$ we get $L_4=\{1324, 1423, 1342, 1432,
3142, 4132\}$, and so on. It turns out that $L_p$ is, in a sense,
an analog of $132$ for $p>3$. For example, the following result is
an analog of Theorem \ref{tcw}~(ii) for  $p=4$.

\begin{theorem}\label{e517}  For any $k\geq2$,
$$
F_{\{L_4,[k]\}}(x)=1+x+x^2R_k(x)R_{k-1}(x)\bigl(R_{k-1}(x)+R_{k-2}(x)\bigr).
$$
\end{theorem}
\begin{proof} The main ingredient of the proof is the following
description of the block decompositions of permutations in
$\SS_n(L_4)$. Let $\alpha\in \SS_n(L_4)$, then there exist $0\ls
r\leq s\leq n-1$ such that either
$\alpha=\alpha_1,n-1,\alpha_2,n,\alpha_3$, or
$\alpha=\alpha_1,n,\alpha_2,n-1,\alpha_3$, where $\alpha_1$ is a
permutation of the numbers $s+1,s+2,\dots,n-2$, $\alpha_2$ is a
permutation of the numbers $r+1,r+2,\dots,s$, and $\alpha_3$ is a
permutation of the numbers $1,2,\dots, r$.
\end{proof}

Consider now the case of a general $p>3$. For an arbitrary
$\pi\in\SS_p$ we define a sequence $a(\pi)$ of zeros and ones of
length $p-1$ as follows. First of all, we put $a_1(\pi)=1$ if
$\pi_{p-1}<\pi_p$ and $a_1(\pi)=0$ otherwise. If
$a_1(\pi),\dots,a_j(\pi)$ are already determined, we put
$a_{j+1}(\pi)=1$ if the length of the maximal increasing
subsequence in $(\pi_{k-j-1},\pi_{k-j},\dots,\pi_k)$ is greater
than that of $(\pi_{k-j},\dots,\pi_k)$, and $a_{j+1}(\pi)=0$
otherwise. For example, let $\pi=7346215$, then $a_1(\pi)=1$,
$a_2(\pi)=0$, $a_3(\pi)=0$, $a_4(\pi)=0$, $a_5(\pi)=1$,
$a_6(\pi)=0$. For an arbitrary sequence $a=(a_1,\dots,a_{p-1})\in
Q^{p-1}=\{0,1\}^{p-1}$ we denote by $N(a)$ the number of
permutations $\pi\in\SS_p$ such that $a(\pi)=a$. The following
result is a further generalization of Theorem \ref{e517}.

\begin{theorem}\label{e518} For any $k\geq p-2$ and any $p>3$,
$$
F_{\{L_p,[k]\}}(x)=\sum_{i=0}^{p-3}i!x^i+x^{p-2}R_k(x)R_{k-1}(x)
\sum_{a\in Q^{p-3}}N(a)\prod_{j=1}^{p-3}R_{k-j-a_{j}}(x).
$$
\end{theorem}

As a generalization of Theorem \ref{a51} we get the following
result.

\begin{theorem}\label{e519} {\rm (i)} For any $m$, $1\leq m\leq k-2$,
$$\begin{array}{l}
F_{\{L_4,[k,m]\}}(x)=1+\frac1{4t^2}+\frac{U_{k-2}(t)U_{m-1}(t)}
{tU_k(t)U_m(t)}\\
\quad\quad\quad+\frac{U_{k-m-2}(t)}{2tU_k(t)U_m(t)}\left(\frac{U_{k-m-2}(t)}{U_{k-m-1}(t)}+
\frac{U_{k-m-3}(t)}{U_{k-m-2}(t)}\right), \qquad
t=\frac1{2\sqrt{x}}.
\end{array}
$$

{\rm (ii)}
$$
F_{\{L_4,[k,k-1]\}}(x)=F_{\{L_4,[k]\}}(x).
$$
\end{theorem}
%=========================================================================
\section*{\bf 6. Counting $321$-restricted permutations}
\setcounter{section}{6} \setcounter{subsection}{0}
\setcounter{theorem}{0}

The case of $321$-restricted permutations is studied much less
than the previous one. We start from the following generalization
of Theorem~\ref{tcw}(i).

\begin{theorem}\label{u1} {\rm (\cite[Th.~1.2(ii)]{MV3})} For $k\geq2$
and any $m$, $1\leq m\leq k-1$,
$$
F_{\{321,[k,m]\}}(x)=R_k(x).
$$
\end{theorem}

Therefore, the Wilf class of $\{132,[k,m]\}$ contains the pair
$\{321,[k,m]\}$ as well. Once again, we have only an analytical
proof of this result. Moreover, the ideas behind this proof are
very different from these behind the proofs in the previous
section. Let $T=\{321, [k,m]\}$; we define
    $$A(n,r)=\sum_{i=0}^{r+m}(-1)^i\binom{r+m-i}{i}f_T(n-i),$$
where $f_T(n)=|S_n(T)|$, and prove that $A(n,r)$ satisfy a linear
recurrence (see \cite[Th.~2.3]{MV3}). It follows that
    $$\sum_{i=0}^k
    (-x)^i\binom{k-i}{i}\left(F_T(x)-\sum_{j=0}^{k-i-1}x^jc_j\right)=0,$$
where $c_j$ is the $j$th Catalan number. Using classical
identities involving Catalan numbers we get the desired result.\\

In view of Theorems \ref{a51} and \ref{a52}, it is a challenge to
find a bijective proof of Theorem \ref{u1}.

A striking analog of Theorem \ref{b57} is given by the following
result.

\begin{theorem} {\rm (\cite[Th.~10]{Kr})}\label{u2} For $k\geq 3$ and
any $r$, $1\leq r \leq k$,
$$
G^r_{321;[k,1]}(x)=\frac{U^{r-1}_{k-1}(t)}
{(2t)^{r-1}U^{r+1}_{k}(t)}, \qquad t=\frac1{2\sqrt{x}}.
$$
\end{theorem}

The proof of this result is based on the bijection $\Psi$ between
$123$-avoiding permutations with exactly $r$ occurrences of the
pattern $[k,1]$ and Dyck paths which start at the origin, return
to the $x$-axis, and have exactly $r$ peaks at height $k$.\\

The case of a general $2$-layered pattern remains intractable. Our
computational experiments suggest that $G^r_{321;
[k,1]}(x)=G^r_{321; [k,2]}(x)$ for $1\ls r \ls k$; however, we are
unable to prove this.
%=========================================================================

\end{document}